\expandafter\chardef\csname pre amssym.def at\endcsname=\the\catcode`\@
\catcode`\@=11

\def\undefine#1{\let#1\undefined}
\def\newsymbol#1#2#3#4#5{\let\next@\relax
 \ifnum#2=\@ne\let\next@\msafam@\else
 \ifnum#2=\tw@\let\next@\msbfam@\fi\fi
 \mathchardef#1="#3\next@#4#5}
\def\mathhexbox@#1#2#3{\relax
 \ifmmode\mathpalette{}{\m@th\mathchar"#1#2#3}%
 \else\leavevmode\hbox{$\m@th\mathchar"#1#2#3$}\fi}
\def\hexnumber@#1{\ifcase#1 0\or 1\or 2\or 3\or 4\or 5\or 6\or 7\or 8\or
 9\or A\or B\or C\or D\or E\or F\fi}

\font\tenmsa=msam10
\font\sevenmsa=msam7
\font\fivemsa=msam5
\newfam\msafam
\textfont\msafam=\tenmsa
\scriptfont\msafam=\sevenmsa
\scriptscriptfont\msafam=\fivemsa
\edef\msafam@{\hexnumber@\msafam}
\mathchardef\dabar@"0\msafam@39
\def\dashrightarrow{\mathrel{\dabar@\dabar@\mathchar"0\msafam@4B}}
\def\dashleftarrow{\mathrel{\mathchar"0\msafam@4C\dabar@\dabar@}}

\def\ulcorner{\delimiter"4\msafam@70\msafam@70 }
\def\urcorner{\delimiter"5\msafam@71\msafam@71 }
\def\llcorner{\delimiter"4\msafam@78\msafam@78 }
\def\lrcorner{\delimiter"5\msafam@79\msafam@79 }
\def\yen{{\mathhexbox@\msafam@55 }}
\def\checkmark{{\mathhexbox@\msafam@58 }}
\def\circledR{{\mathhexbox@\msafam@72 }}
\def\maltese{{\mathhexbox@\msafam@7A }}

\font\tenmsb=msbm10
\font\sevenmsb=msbm7
\font\fivemsb=msbm5
\newfam\msbfam
\textfont\msbfam=\tenmsb
\scriptfont\msbfam=\sevenmsb
\scriptscriptfont\msbfam=\fivemsb
\edef\msbfam@{\hexnumber@\msbfam}

\catcode`\@=\csname pre amssym.def at\endcsname

\expandafter\ifx\csname pre amssym.tex at\endcsname\relax \else \endinput\fi
\expandafter\chardef\csname pre amssym.tex at\endcsname=\the\catcode`\@
\catcode`\@=11
\newsymbol\boxdot 1200
\newsymbol\boxplus 1201
\newsymbol\boxtimes 1202
\newsymbol\square 1003
\newsymbol\blacksquare 1004
\newsymbol\centerdot 1205
\newsymbol\lozenge 1006
\newsymbol\blacklozenge 1007
\newsymbol\circlearrowright 1308
\newsymbol\circlearrowleft 1309
\undefine\rightleftharpoons
\newsymbol\rightleftharpoons 130A
\newsymbol\leftrightharpoons 130B
\newsymbol\boxminus 120C
\newsymbol\Vdash 130D
\newsymbol\Vvdash 130E
\newsymbol\vDash 130F
\newsymbol\twoheadrightarrow 1310
\newsymbol\twoheadleftarrow 1311
\newsymbol\leftleftarrows 1312
\newsymbol\rightrightarrows 1313
\newsymbol\upuparrows 1314
\newsymbol\downdownarrows 1315
\newsymbol\upharpoonright 1316
 
\newsymbol\downharpoonright 1317
\newsymbol\upharpoonleft 1318
\newsymbol\downharpoonleft 1319
\newsymbol\rightarrowtail 131A
\newsymbol\leftarrowtail 131B
\newsymbol\leftrightarrows 131C
\newsymbol\rightleftarrows 131D
\newsymbol\Lsh 131E
\newsymbol\Rsh 131F
\newsymbol\rightsquigarrow 1320
\newsymbol\leftrightsquigarrow 1321
\newsymbol\looparrowleft 1322
\newsymbol\looparrowright 1323
\newsymbol\circeq 1324
\newsymbol\succsim 1325
\newsymbol\gtrsim 1326
\newsymbol\gtrapprox 1327
\newsymbol\multimap 1328
\newsymbol\therefore 1329
\newsymbol\because 132A
\newsymbol\doteqdot 132B
 
\newsymbol\triangleq 132C
\newsymbol\precsim 132D
\newsymbol\lesssim 132E
\newsymbol\lessapprox 132F
\newsymbol\eqslantless 1330
\newsymbol\eqslantgtr 1331
\newsymbol\curlyeqprec 1332
\newsymbol\curlyeqsucc 1333
\newsymbol\preccurlyeq 1334
\newsymbol\leqq 1335
\newsymbol\leqslant 1336
\newsymbol\lessgtr 1337
\newsymbol\backprime 1038
\newsymbol\risingdotseq 133A
\newsymbol\fallingdotseq 133B
\newsymbol\succcurlyeq 133C
\newsymbol\geqq 133D
\newsymbol\geqslant 133E
\newsymbol\gtrless 133F
\newsymbol\sqsubset 1340
\newsymbol\sqsupset 1341
\newsymbol\vartriangleright 1342
\newsymbol\vartriangleleft 1343
\newsymbol\trianglerighteq 1344
\newsymbol\trianglelefteq 1345
\newsymbol\bigstar 1046
\newsymbol\between 1347
\newsymbol\blacktriangledown 1048
\newsymbol\blacktriangleright 1349
\newsymbol\blacktriangleleft 134A
\newsymbol\vartriangle 134D
\newsymbol\blacktriangle 104E
\newsymbol\triangledown 104F
\newsymbol\eqcirc 1350
\newsymbol\lesseqgtr 1351
\newsymbol\gtreqless 1352
\newsymbol\lesseqqgtr 1353
\newsymbol\gtreqqless 1354
\newsymbol\Rrightarrow 1356
\newsymbol\Lleftarrow 1357
\newsymbol\veebar 1259
\newsymbol\barwedge 125A
\newsymbol\doublebarwedge 125B
\undefine\angle
\newsymbol\angle 105C
\newsymbol\measuredangle 105D
\newsymbol\sphericalangle 105E
\newsymbol\varpropto 135F
\newsymbol\smallsmile 1360
\newsymbol\smallfrown 1361
\newsymbol\Subset 1362
\newsymbol\Supset 1363
\newsymbol\Cup 1264
 
\newsymbol\Cap 1265
 
\newsymbol\curlywedge 1266
\newsymbol\curlyvee 1267
\newsymbol\leftthreetimes 1268
\newsymbol\rightthreetimes 1269
\newsymbol\subseteqq 136A
\newsymbol\supseteqq 136B
\newsymbol\bumpeq 136C
\newsymbol\Bumpeq 136D
\newsymbol\lll 136E
 
\newsymbol\ggg 136F
 
\newsymbol\circledS 1073
\newsymbol\pitchfork 1374
\newsymbol\dotplus 1275
\newsymbol\backsim 1376
\newsymbol\backsimeq 1377
\newsymbol\complement 107B
\newsymbol\intercal 127C
\newsymbol\circledcirc 127D
\newsymbol\circledast 127E
\newsymbol\circleddash 127F
\newsymbol\lvertneqq 2300
\newsymbol\gvertneqq 2301
\newsymbol\nleq 2302
\newsymbol\ngeq 2303
\newsymbol\nless 2304
\newsymbol\ngtr 2305
\newsymbol\nprec 2306
\newsymbol\nsucc 2307
\newsymbol\lneqq 2308
\newsymbol\gneqq 2309
\newsymbol\nleqslant 230A
\newsymbol\ngeqslant 230B
\newsymbol\lneq 230C
\newsymbol\gneq 230D
\newsymbol\npreceq 230E
\newsymbol\nsucceq 230F
\newsymbol\precnsim 2310
\newsymbol\succnsim 2311
\newsymbol\lnsim 2312
\newsymbol\gnsim 2313
\newsymbol\nleqq 2314
\newsymbol\ngeqq 2315
\newsymbol\precneqq 2316
\newsymbol\succneqq 2317
\newsymbol\precnapprox 2318
\newsymbol\succnapprox 2319
\newsymbol\lnapprox 231A
\newsymbol\gnapprox 231B
\newsymbol\nsim 231C
\newsymbol\ncong 231D
\newsymbol\diagup 231E
\newsymbol\diagdown 231F
\newsymbol\varsubsetneq 2320
\newsymbol\varsupsetneq 2321
\newsymbol\nsubseteqq 2322
\newsymbol\nsupseteqq 2323
\newsymbol\subsetneqq 2324
\newsymbol\supsetneqq 2325
\newsymbol\varsubsetneqq 2326
\newsymbol\varsupsetneqq 2327
\newsymbol\subsetneq 2328
\newsymbol\supsetneq 2329
\newsymbol\nsubseteq 232A
\newsymbol\nsupseteq 232B
\newsymbol\nparallel 232C
\newsymbol\nmid 232D
\newsymbol\nshortmid 232E
\newsymbol\nshortparallel 232F
\newsymbol\nvdash 2330
\newsymbol\nVdash 2331
\newsymbol\nvDash 2332
\newsymbol\nVDash 2333
\newsymbol\ntrianglerighteq 2334
\newsymbol\ntrianglelefteq 2335
\newsymbol\ntriangleleft 2336
\newsymbol\ntriangleright 2337
\newsymbol\nleftarrow 2338
\newsymbol\nrightarrow 2339
\newsymbol\nLeftarrow 233A
\newsymbol\nRightarrow 233B
\newsymbol\nLeftrightarrow 233C
\newsymbol\nleftrightarrow 233D
\newsymbol\divideontimes 223E
\newsymbol\varnothing 203F
\newsymbol\nexists 2040
\newsymbol\Finv 2060
\newsymbol\Game 2061
\newsymbol\mho 2066
\newsymbol\eth 2067
\newsymbol\eqsim 2368
\newsymbol\beth 2069
\newsymbol\gimel 206A
\newsymbol\daleth 206B
\newsymbol\lessdot 236C
\newsymbol\gtrdot 236D
\newsymbol\ltimes 226E
\newsymbol\rtimes 226F
\newsymbol\shortmid 2370
\newsymbol\shortparallel 2371
\newsymbol\smallsetminus 2272
\newsymbol\thicksim 2373
\newsymbol\thickapprox 2374
\newsymbol\approxeq 2375
\newsymbol\succapprox 2376
\newsymbol\precapprox 2377
\newsymbol\curvearrowleft 2378
\newsymbol\curvearrowright 2379
\newsymbol\digamma 207A
\newsymbol\varkappa 207B
\newsymbol\Bbbk 207C
\newsymbol\hslash 207D
\undefine\hbar
\newsymbol\hbar 207E
\newsymbol\backepsilon 237F
\catcode`\@=\csname pre amssym.tex at\endcsname

\magnification=1200
\hsize=468truept
\vsize=646truept
\voffset=-10pt
\parskip=4pt
\baselineskip=14truept
\count0=1

\dimen100=\hsize

\def\leftill#1#2#3#4{
\medskip
\line{$
\vcenter{
\hsize = #1truept \hrule\hbox{\vrule\hbox to  \hsize{\hss \vbox{\vskip#2truept
\hbox{{\copy100 \the\count105}: #3}\vskip2truept}\hss }
\vrule}\hrule}
\dimen110=\dimen100
\advance\dimen110 by -36truept
\advance\dimen110 by -#1truept
\hss \vcenter{\hsize = \dimen110
\medskip
\noindent { #4\par\medskip}}$}
\advance\count105 by 1
}
\def\rightill#1#2#3#4{
\medskip
\line{
\dimen110=\dimen100
\advance\dimen110 by -36truept
\advance\dimen110 by -#1truept
$\vcenter{\hsize = \dimen110
\medskip
\noindent { #4\par\medskip}}
\hss \vcenter{
\hsize = #1truept \hrule\hbox{\vrule\hbox to  \hsize{\hss \vbox{\vskip#2truept
\hbox{{\copy100 \the\count105}: #3}\vskip2truept}\hss }
\vrule}\hrule}
$}
\advance\count105 by 1
}
\def\midill#1#2#3{\medskip
\line{$\hss
\vcenter{
\hsize = #1truept \hrule\hbox{\vrule\hbox to  \hsize{\hss \vbox{\vskip#2truept
\hbox{{\copy100 \the\count105}: #3}\vskip2truept}\hss }
\vrule}\hrule}
\dimen110=\dimen100
\advance\dimen110 by -36truept
\advance\dimen110 by -#1truept
\hss $}
\advance\count105 by 1
}
\def\insectnum{\copy110\the\count120
\advance\count120 by 1
}

\font\ninerm=cmr9
\font\eightrm=cmr8

\font\tenrm=cmr10 at 10pt

\font\sc=cmcsc10

\def\msb{\fam\msbfam\tenmsb}

\def\bbc{{\msb C}}

\def\bbi{{\msb I}}

\def\bbp{{\msb P}}
\def\bbq{{\msb Q}}
\def\bbr{{\msb R}}

\def\bbz{{\msb Z}}

\def\grD{\Delta}

\def\grL{\Lambda}

\def\grS{\Sigma}

\def\gra{\alpha}

\def\grl{\lambda}

\def\gro{\omega}

\def\grs{\sigma}

\def\la#1{\hbox to #1pc{\leftarrowfill}}
\def\ra#1{\hbox to #1pc{\rightarrowfill}}

\def\fract#1#2{\raise4pt\hbox{$ #1 \atop #2 $}}
\def\decdnar#1{\phantom{\hbox{$\scriptstyle{#1}$}}
\left\downarrow\vbox{\vskip15pt\hbox{$\scriptstyle{#1}$}}\right.}

\def\bowtie{\hbox to 1pt{\hss}\raise.66pt\hbox{$\scriptstyle{>}$}
\kern-4.9pt\triangleleft}
\def\hsmash{\triangleright\kern-4.4pt\raise.66pt\hbox{$\scriptstyle{<}$}}
\def\boxit#1{\vbox{\hrule\hbox{\vrule\kern3pt
\vbox{\kern3pt#1\kern3pt}\kern3pt\vrule}\hrule}}

\def\za{\vrule height6pt width4pt depth1pt}

\font\aa=eufm10

\def\Got#1{\hbox{\aa#1}}

\def\bfw{{\bf w}}

\def\calo{{\cal O}}

\def\calf{{\cal F}}

\def\calo{{\cal O}}

\def\cals{{\cal S}}

\def\calz{{\cal Z}}

\def\gg{{\Got g}}

\def\gp{{\Got p}}

\def\Got#1{\hbox{\aa#1}}

\def\gsp1{{\Got s}{\Got p}(1)}

\font\svtnrm=cmr17

\font\bsc=cmcsc10 at 10truept

\def\coker{\hbox{coker}}
\def\ker{\hbox{ker}}

\def\Se{Sasakian-Einstein }

\def\mult{\hbox{mult}}

\def\link{2}

\def\top{3}
\def\dis{4}
\def\reg{5}

\phantom{ooo}
\bigskip\bigskip
\centerline{\svtnrm Einstein Metrics on} \medskip
\centerline{\svtnrm Rational Homology 7-Spheres}  \medskip

\bigskip\bigskip
\centerline{\sc Charles P. Boyer~~ Krzysztof Galicki~~ Michael Nakamaye}
\footnote{}{\ninerm During the preparation of this work the first two authors 
were partially supported by NSF grant DMS-9970904, and third author by NSF
grant DMS-0070190.} \bigskip

\centerline{\vbox{\hsize = 5.85truein
\baselineskip = 12.5truept
\eightrm
\noindent {\bsc Abstract:} In this paper we demonstrate the existence of \Se
structures on certain 2-connected rational homology 7-spheres.  These appear 
to be the first non-regular examples of  \Se metrics on simply connected 
rational homology spheres. We also briefly describe the rational homology 
7-spheres that admit regular positive Sasakian structures. }}\tenrm

\vskip .5in
\baselineskip = 10 truept
\centerline{\bf Introduction}  
\bigskip

Dimension seven appears to be rather special when it comes to examples of
compact Einstein manifolds. It is perhaps the
prominent r\^ole such manifolds have played in physics ever since the early
days of Kaluza-Klein supergravity that made both theoretical physicists and
mathematicians alike particularly interested in them. Discoveries of many
different constructions followed as a result of this interest.

Arguably, today a 
special place among all compact Einstein 7-manifolds is reserved for
the so-called Sasakian-Einstein spaces. They are defined to be
Riemannian manifolds with the property that the metric cone on them is
a Calabi-Yau 4-fold and, in particular, they are always of
positive scalar curvature. All regular \Se manifolds are circle bundles of
Fano 3-folds that admit K\"ahler-Einstein metrics. Non-regular ones
fiber over compact K\"ahler-Einstein Fano 3-folds with orbifold singularities.
An interesting sub-family of the family of \Se 7-manifolds consists of
the so-called 3-Sasakian spaces. They are characterized by fact that
their metric cone is 
not only Calabi-Yau, but also hyperk\"ahler and are all
orbifold fibrations over compact K\"ahler-Einstein Fano 3-folds which admit a
complex contact structure.

Regular and non-regular examples of both \Se and 3-Sasakian manifolds are now
plentiful and they were extensively studied by the first two authors [BG1,
BG2]. There is an example of a regular \Se (4n+3)-manifold
which is worthy of some further discussion. It is the homogeneous Stiefel
manifold of 2-frames in (2n+1)-dimensional Euclidean space,
$V_2(\bbr^{2n+1})=SO(2n+1)/SO(2n-1)$ which is a circle bundle over the oriented
Grassmannian, $\widetilde{G}_2(\bbr^{2n+1}).$ From the point of view of an
algebraic geometer it is a classical fact that $\widetilde{G}_2(\bbr^{2n+1})$
is diffeomorphic to  the complex quadric $Q_{2n-1}$ in $\bbc\bbp^{2n}$ which
is well-known to be Fano and to admit a K\"ahler-Einstein metric. It is
perhaps less well-known that the quadric $Q_{2n-1}$ has the same cohomology
groups as $\bbc\bbp^{2n-1}$ , but differs in the ring structure. Hence,
$V_2(\bbr^{2n+1})$ is a rational  homology sphere with
$H_{n-1}(V_2(\bbr^{2n+1}),\bbz)\approx \bbz_2.$ Now it has been known for quite
some time that $V_2(\bbr^{2n+1})$ carries a \Se structure [BGFK,BG1]. Up to
date, apart from $S^{2n+1}, V_2(\bbr^{2n+1})$, and the 3-Sasakian homogeneous
11-manifold $G_2/Sp(1)$ we are not aware of any other examples of simply
connected rational homology spheres which are also known to admit \Se
structures. In this paper we shall demonstrate that for $n=2,$ quite to the
contrary, there are many examples of such structures, 184 to be precise. These
examples are obtained as hypersurfaces in certain weighted projective
4-spaces, but we certainly expect the phenomena to occur in arbitrary
dimension.

The key to this construction is a recent paper of Johnson and Koll\'ar [JK2].
There they give a list of 4442 quasi-smooth Fano
3-folds $\calz$ anticanonically embedded in weighted projective 4-spaces
$\bbp(\bfw).$  Moreover, they show that 1936 of these 3-folds admit
K\"ahler-Einstein metrics. According to our general theory [BG1] such Fano
3-folds give rise to \Se metrics on smooth 7-manifolds $M^7.$ Moreover, these
7-manifolds arise as links of isolated hypersurface singularities associated to
certain weighted homogeneous polynomials in $\bbc^5.$ As in [JK1] Johnson and
Koll\'ar [JK2] only consider the case when the orbifold Fano index is one, and
as the authors showed in [BGN1] for log del Pezzo surfaces, there should be
many more interesting examples of quasi-smooth Fano 3-folds with higher
orbifold Fano index. This is currently under study. 

In this note we prove

\noindent{\sc Theorem A}: \tensl There are1936 distinct \Se structures on 
certain $2$-connected 
\vskip1pt \noindent 7-manifolds $M^7_{\bfw,d}$ realized as links of weighted 
homogeneous polynomials in $\bbc^5$ with weight vector 
$\bfw=(w_0,w_1,w_2,w_3,w_4)$ and degree $d.$ In particular, there are 184  
$2$-connected rational homology spheres which are listed in the Table below. 
In addition to the weight vector $\bfw$ and degree $d$ the Table lists the 
Milnor number $\mu$ of the link and the order of $H_3(M^7_{\bfw,d},\bbz).$
\tenrm

We have not answered the question as to whether two distinct or non-conjugate 
\Se structures on the same link $M^7$ could belong to the same underlying 
Riemannian metric $g.$ Indeed, this can happen, but if $g$ is not the standard 
round metric on $S^7$ then by a Theorem of Tachibana and Yu [TaYu], the two 
\Se structures must belong to a 3-Sasakian structure. But then by a Theorem of 
Galicki and Salamon [GS],  we must have $b_3=0,$ so $M^7$ must be a rational 
homology sphere. However, we do not know whether any of the rational homology 
7-spheres discussed here admit 3-Sasakian structures.

\settabs 10\columns
\bigskip
\bigskip
\baselineskip = 10 truept
\centerline{\bf \link. The Sasakian Geometry of Links of Weighted
Homogeneous Polynomials}   \bigskip

In this section we briefly review the Sasakian geometry of links of isolated
hypersurface singularities defined by weighted homogeneous polynomials.
Consider the affine space $\bbc^{n+1}$ together with a weighted
$\bbc^*$-action given by $(z_0,\ldots,z_n)\mapsto
(\grl^{w_0}z_0,\ldots,\grl^{w_n}z_n),$ where the {\it weights} $w_j$ are
positive integers. It is convenient to view the weights as the components of a
vector $\bfw\in (\bbz^+)^{n+1},$ and we shall assume that they are ordered
$w_0\leq w_1\leq \cdots\leq w_n$ and that $\gcd(w_0,\ldots,w_n)=1.$ Let $f$
be a quasi-homogeneous polynomial, that is $f\in \bbc[z_0,\ldots,z_n]$ and
satisfies $$f(\grl^{w_0}z_0,\ldots,\grl^{w_n}z_n)=\grl^df(z_0,\ldots,z_n),
\leqno{\link.1}$$
where $d\in \bbz^+$ is the degree of $f.$ We are interested in the {\it
weighted affine cone} $C_f$ defined by
the equation $f(z_0,\ldots,z_n)=0.$ We shall assume that the origin in
$\bbc^{n+1}$ is an isolated singularity, in fact the only singularity, of
$f.$ Then the link $L_f$ defined by 
$$L_f= C_f\cap S^{2n+1}, \leqno{\link.2}$$
where 
$$S^{2n+1}=\{(z_0,\ldots,z_n)\in \bbc^{n+1}|\sum_{j=0}^n|z_j|^2=1\}$$
is the unit sphere in $\bbc^{n+1},$ is a smooth manifold of dimension $2n-1.$ 
Furthermore, it is well-known [Mil] that the link $L_f$ is $(n-2)$-connected.

On $S^{2n+1}$ there is a well-known [YK] ``weighted'' Sasakian structure  
$(\xi_\bfw,\eta_\bfw,\Phi_\bfw,g_\bfw)$ which in the standard coordinates
$\{z_j=x_j+iy_j\}_{j=0}^n$ on $\bbc^{n+1}=\bbr^{2n+2}$ is determined by
$$\eta_\bfw = {\sum_{i=0}^n(x_idy_i-y_idx_i)\over\sum_{i=0}^n
w_i(x_i^2+y_i^2)}, \qquad \xi_\bfw
=\sum_{i=0}^nw_i(x_i\partial_{y_i}-y_i\partial_{x_i}),$$
and the standard Sasakian structure $(\xi,\eta,\Phi,g)$ on $S^{2n+1}.$
The embedding $L_f\hookrightarrow S^{2n+1}$ induces a Sasakian structure on
$L_f$ [BG3]. 

Given a sequence $\bfw =(w_0,\ldots,w_n)$ of ordered positive integers one can
form the graded polynomial ring $S(\bfw)=\bbc[z_0,\ldots,z_n]$, where $z_i$ has
grading or {\it weight} $w_i.$ The weighted projective space [Dol, Fle]
$\bbp(\bfw)=\bbp(w_0,\ldots,w_n)$ is defined to be the scheme
$\hbox{Proj}(S(\bfw)).$  It is the quotient space 
$(\bbc^{n+1}-\{0\})/\bbc^*(\bfw)$, where $\bbc^*(\bfw)$ is the weighted action
defined in \link.1, or equivalently, $\bbp(\bfw)$ is the quotient of
the weighted Sasakian sphere
$S_\bfw^{2n+1}=(S^{2n+1},\xi_\bfw,\eta_\bfw,\Phi_\bfw,g_\bfw)$ by the
weighted circle action $S^1(\bfw)$ generated by $\xi_\bfw.$ As such
$\bbp(\bfw)$ is also a compact complex orbifold with an induced K\"ahler
structure. We have from [BG3]

\noindent{\sc Theorem} \link.3: \tensl The quadruple
$(\xi_\bfw,\eta_\bfw,\Phi_\bfw,g_\bfw)$ gives $L_f$ a quasi-regular Sasakian
structure such that there is a commutative diagram
$$\matrix{L_f &\ra{2.5}& S^{2n+1}_\bfw&\cr
  \decdnar{\pi}&&\decdnar{} &\cr
   \calz_f &\ra{2.5} &\bbp(\bfw),&\cr}$$
where the horizontal arrows are Sasakian and K\"ahlerian embeddings,
respectively, and the vertical arrows are principal $S^1$ V-bundles and
orbifold Riemannian submersions.   Moreover, if $\calz_f$ is Fano, $L_f$ is 
the total space of the principal $S^1$ V-bundle over the orbifold $\calz_f$ 
whose first Chern class  in $H^2_{orb}(\calz_f,\bbz)$ is $c_1(\calz_f)/I,$ 
where $I$ is the index. \tenrm 

We should also mention that 
$c_1(\calz_f)$ pulls back to the basic first Chern
class $c_1^B\in H^2_B(\calf_{\xi_\bfw})$ and $\eta_\bfw$ is the connection in
this V-bundle whose curvature is $d\eta ={2ni\over I}\pi^*\gro_\bfw,$ where
$\gro_\bfw$ is the K\"ahler form on $\calz_f.$ 

Now conditions on the weights that guarantee that the hypersurface $C_f\subset
\bbc^{n+1}$ have only an isolated singularity at the origin are well-known
[Fle,JK1]. These conditions become more complicated as the dimension increases
[Fle,JK2]; however, in this paper we are only interested in the $n=4$ case
of hypersurfaces in a weighted complex projective 4-space.  These conditions,
known as {\it quasi-smoothness} conditions guarantee that $\calz_f$ is smooth
in the orbifold sense, that is, at a vertex  $P_i\in\bbp({\bf w})$ the preimage
of $\calz_f$ in the orbifold chart of  $\bbp({\bf w})$ is smooth.  It is
easy to see that one can formulate all these conditions as follows [Fle,JK2]: 
\bigskip
\noindent{\sc Quasi-Smoothness Conditions} \link.4:
\bigskip
\+I.&For each $i=0,\cdots,4$ there is a $j$ and a monomial
$z_i^{m_i}z_j\in \calo(d).$\cr
\+&Here $j=i$ is possible.\cr
\medskip
\+II.&For all distinct $i,j$ either there is a 
monomial $z_i^{b_i}z_j^{b_j}\in
\calo(d).$\cr
\+&or there exist monomials $z^{n_1}_iz^{m_1}_jz_k,z^{n_2}_iz^{m_2}_jz_l
\in \calo(d)$ with $\{k,l\}\neq \{i,j\}$ and $k\neq l.$\cr
\medskip
\+III.&For every $i,j$  there exists a monomial of degree $d$ that does not
involve either $z_i$\cr 
\+&or $z_j.$ \cr

There is another condition apart from quasi-smoothness that assures us
that the adjunction theory behaves correctly, and that $\bbp({\bf w})$ 
does not have any orbifold singularities of codimension 1. It is [Dol,Fle] 
\medskip
\noindent{\sc Well-formedness Condition} \link.5
\bigskip
\+IV.& For each $i$ we have
$\gcd(w_0,\cdots,\hat{w_i},\cdots,w_4)=1.$ \cr
\+&Here the $\hat{}$ means skip that element.\cr \medskip

Condition IV
guarantees that the canonical V-bundle $K_\calz$ is
determined in terms of the degree and index by
$$K_\calz \simeq \calo(-I)=\calo(d-|\bfw|),\leqno{\link.6}$$
where $|\bfw|=\sum_iw_i.$ 

In this note we shall only consider the anticanonically embedded Fano 3-folds
of [JK2], that is, we shall assume hereafter that $I=|\bfw|-d=1.$ The examples
we consider are from the list sporadic.txt of Johnson and Koll\'ar [JK2] which
is found at: 

\noindent http://www.math.princeton.edu/\~~jmjohnso. 

\bigskip
\centerline{\bf \top. The Topology of the Link $M^7_{\bfw,d}$}   
\bigskip

The topology of a link $L_f$ of an isolated hypersurface singularity is
encoded in the  characteristic polynomial $\grD(t)$ of the monodromy map.
$\grD(t)$ is an important link invariant that generalizes the Alexander
polynomial of a knot, and is often called the ``Alexander polynomial'' of the
link [HZ].  Let us  recall the well-known construction of Milnor [Mil]
concerning isolated hypersurface singularities: There is a fibration of
$(S^{2n+1}-L_f)\ra{1.3} S^1$ whose fiber $F$ is an open manifold that is
homotopy equivalent to a bouquet of n-spheres $S^n\vee S^n\cdots \vee S^n.$
The {\it Milnor number} $\mu$ of $L_f$ is the number of $S^n$'s in the
bouquet. It is an invariant of the link which can be calculated explicitly in
terms of the degree $d$ and weights $(w_0,\ldots,w_n)$ by the formula [MO]    
$$\mu =\mu(L_f)=\prod_{i=0}^n\bigl({d\over w_i}-1\bigr).\leqno{\top.1}$$   
The closure $\bar{F}$ of $F$ has the same
homotopy type as $F$ and is a compact manifold whose boundary is precisely the
link $L_f.$ So the reduced homology of $F$ and $\bar{F}$ is only non-zero in
dimension $n$ and $H_n(F,\bbz)\approx \bbz^{\mu}.$ Using the Wang sequence of
the Milnor fibration together with Alexander-Poincare duality gives the exact
sequence [Mil] $$0\ra{1.5} H_n(L_f,\bbz)\ra{1.5} H_n(F,\bbz)
\fract{\bbi -h_*}{\ra{1.5}} H_n(F,\bbz) \ra{1.5} H_{n-1}(L_f,\bbz)\ra{1.5} 0,
\leqno{\top.2}$$ 
where $h_*$ is the {\it monodromy} map (or characteristic
map) induced by the $S^1_\bfw$ action. From this we see that
$H_n(L_f,\bbz)=\ker(\bbi -h_*)$ is a free Abelian group, and $H_{n-1}(L_f,\bbz)
=\coker(\bbi -h_*)$ which in general has torsion, but whose free part equals 
$\ker(\bbi -h_*).$ There is a well-known algorithm due to Milnor and Orlik [MO]
for computing the free part of $H_{n-1}(L_f,\bbz)$ in terms of the
characteristic polynomial $\grD(t)=\det(t\bbi -h_*)$ of the monodromy map. The
Betti number  $b_n(L_f)=b_{n-1}(L_f)$ equals the number of factors of $(t-1)$
in $\grD(t).$   Generally, finding the torsion is much more difficult. However,
in the case of rational homology spheres, $b_n(L_f)=b_{n-1}(L_f)=0,$ the group
$H_{n-1}(M,\bbz)$ is a torsion group of order $\grD(1).$

It is not our purpose in this note to give a systematic study of the
Johnson-Koll\'ar list. This requires a computer program for computing the
Betti numbers which is currently under study. Here we are content with giving
algorithm for finding special cases when rational homology spheres occur. We
have written a MAPLE  program which allows us to search the JK list,
sporatic.txt and pick out certain rational homology spheres. We emphasize that
this procedure does not necessarily find all rational homology spheres, but
only all of those that satisfy certain additional conditions. Actually there
are two distinct types of conditions on the weights that allow us to find
rational homology spheres and they are described in the lemmas below. The
first and simplest is that the weights are all relatively prime to the degree.

\noindent{Lemma} \top.4: \tensl Let $\bfw=(w_0,w_1,w_2,w_3,w_4)$ be the
weights of a quasi-smooth Fano 3-fold, $\calz_f$ of degree $d$ and index $1.$
Suppose further that $\gcd(w_i,d)=1$ for all $i=0,\cdots, 4.$ Then there
exists an integer $N(\bfw)$ such that the Alexander polynomial $\grD(t)$
of the link $M^7_{\bfw,d}$ has the form
$$\grD(t)={(t^d-1)^{N(\bfw)}\over t-1}$$
Hence, the Betti number $b_3$ of
the link $M^7_{\bfw,d}$ is given by
$b_3(M^7_{\bfw,d}) = N(\bfw)-1.$
\tenrm

\noindent{\sc Proof}: The Milnor and Orlik [MO] algorithm for computing the
characteristic polynomial of the monodromy operator for weighted homogeneous
polynomials is as follows: First associate to any monic polynomial $F$ with
roots $\gra_1,\ldots,\gra_k\in \bbc^*$ its divisor 
$$\hbox{div}~F= <\gra_1>+\cdots+<\gra_k>$$
as an element of the integral ring $\bbz[\bbc^*]$ and let $\grL_n= \hbox{div}~
(t^n-1).$  The rational weights $w'_i$ used in [MO] are related to our integer
weights $w_i$ by $w_i'={d\over w_i},$ and we write the $w'_i={u_i\over v_i}$ in
irreducible form. The divisor $\hbox{div}~\grD$ is given by
$$\hbox{div}~\grD= \Bigl({\grL_{u_0}\over v_0}-1\Bigr)\cdots
\Bigl({\grL_{u_4}\over v_4}-1\Bigr)\leqno{\top.6}$$
which can be reduced to the form
$$\hbox{div}~ \grD(t)=\sum_ja_j\grL_j-1\leqno{\top.7}$$
for some integers $a_j$ upon using the relations
$\grL_a\grL_b=\gcd(a,b)\grL_{lcm(a,b)}.$  The characteristic polynomial
$\grD(t)$ is then determined from its divisor by
$$\grD(t)={\prod(t^j-1)^{a_j}\over t-1}, \leqno{\top.8}$$
and the third Betti number is given by 
$$b_3(M^7_{\bfw,d})=\sum_ja_j-1.$$

In our case we have $\gcd(w_i,d)=1$ so equation \top.7 must take the form
$$\hbox{div}~ \grD(t)=N(\bfw)\grL_d-1\leqno{\top.9}$$ where $N(\bfw)$ is an
integer. \hfill\za

The integer $N(\bfw)$ can be computed explicitly from the above procedure in
terms of the weights and index. We find 
$$N(\bfw)={d(dr_{01}r_{23}+r_{01}+r_{23})\over w_4}+{1\over w_4}
-(dr_{01}r_{23}+r_{01}+r_{23})\leqno{\top.10}$$
where
$$d=|\bfw|-1,\qquad r_{ij}={d\over w_iw_j}-{1\over w_i}-{1\over w_j}.$$

We should remark here that although it is far from manifest in equation
\top.10, under the hypothesis of Lemma \top.4 the function $N(\bfw)$ is
invariant under a permutation of the weights, i.e. if $\grS_5$ denotes the
permutation group on 5 letters, then $N(\grs(\bfw))=N(\bfw)$ for any $\grs\in
\grS_5.$  We have an immediate 

\noindent{\sc Corollary} \top.11: \tensl Let $M^7_{\bfw,d}$ be the link of an 
isolated
hypersurface singularity defined by a weighted homogeneous polynomial $f$ with
well-formed weights $\bfw=(w_0,w_1,w_2,w_3,w_4)$ and degree $d$ which satisfy
$\gcd(w_i,d)=1$ for all $i=0,\cdots,4.$ Then $M^7_{\bfw,d}$ is a rational 
homology
sphere if and only if $N(\bfw)=1.$ Furthermore, in this case the Milnor number
$\mu=d-1,$ and the order of $H_3(M^7_{\bfw,d},\bbz)$ equals the degree $d.$
\tenrm

\noindent{\sc Proof}: The only part that we need to compute is the order of
$H_3.$ Since for a 2-connected rational homology sphere
$H_4(M^7_{\bfw,d},\bbz)=0,$ the exact sequence \top.2 shows [Mil] that  the
order of $H_3(M^7_{\bfw,d},\bbz)$ equals $\grD(1).$ But from \top.8 and \top.9
we see that in our special case the characteristic polynomial takes the form
$$\grD(t)= {t^d-1\over t-1}=t^{d-1}+\cdots +t+1$$ from which the result
follows. \hfill\za

We now describe the second type of condition.

\noindent{\sc Lemma} \top.12: \tensl Let $\bfw=(w_0,w_1,w_2,w_3,w_4)$ be the
weights of a quasi-smooth Fano 3-fold, $\calz_f$ of degree $d$ and index $1.$
Suppose further that the degree can be written as $d=m_3m_2,$ where $m_2$ and
$m_3$ are relatively prime, and that the ``rational weights'' ${d\over w_i}$
take the form ${m_3\over v_i}$ for 3 values of $i$ and ${m_2\over v_i}$ for the
other 2 values of $i.$ Then there exist positive integers $l=l(\bfw),$
and $n=n(\bfw),$ depending on the weights $\bfw,$ such that the Alexander
polynomial $\grD(t)$ of the link $M^7_{\bfw,d}$ takes the form 
$$\grD(t)={(t^d-1)^{ln}(t^{m_3}-1)^l\over (t-1)(t^{m_2}-1)^n}.$$
Hence, 
$$b_3(M^7_{\bfw,d})=(n(\bfw)+1)(l(\bfw)-1).$$
\tenrm

\noindent{\sc Proof}: Computing as in the proof of Lemma \top.4, we see that
from the Milnor-Orlik procedure [MO] that the divisor of the Alexander
polynomial must take the form
$$\hbox{div}~
\grD(t)=l(\bfw)n(\bfw)\grL_d+l(\bfw)\grL_{m_3}-n(\bfw)\grL_{m_2}-1\leqno{\top.13}$$
for some positive integers $l(\bfw)$ and $n(\bfw)$ depending on the weights.
The above form of the Alexander polynomial then follows from equations \top.7
and \top.8. The explicit form of the functions $l(\bfw)$ and $n(\bfw)$ are
also easily calculated. Let $i_1,i_2,i_3$ denote the 3 indices whose rational
weights take the form ${m_3\over v_i}$ and similarly let $j_1,j_2$ denote the
indices corresponding to the rational weights ${m_2\over v_j}.$ Then one finds
$$\leqalignno{l(\bfw)&={m_3^2\over v_{i_1}v_{i_2}v_{i_3}}-m_3\Bigl({1\over
v_{i_1}v_{i_2}} +{1\over v_{i_1}v_{i_3}} +{1\over v_{i_2}v_{i_3}}\Bigr)
+{1\over v_{i_1}}+ {1\over v_{i_2}} +{1\over v_{i_3}}&\top.14\cr
n(\bfw)&={m_2\over v_{j_1}v_{j_2}}-{1\over v_{j_1}}- {1\over
v_{j_2}}&\top.15\cr}$$
The expression for $b_3$ follows directly from the expression for $\grD(t).$
\hfill\za

\noindent{\sc Corollary} \top.16: \tensl Let $\bfw=(w_0,w_1,w_2,w_3,w_4)$ be
the weights of a quasi-smooth Fano 3-fold, $\calz_f$ of degree $d$ and index
$1.$ Suppose further that the hypothesis of Lemma \top.12 is satisfied, then
$M^7_{\bfw,d}$ is a rational homology sphere if and only if $l(\bfw)=1.$
Furthermore, in this case the Milnor number $\mu=(m_3-1)(n(\bfw)m_2+1),$ and
the order of $H_3(M^7_{\bfw,d},\bbz)$ is $m_3^{n(\bfw)+1}.$ \tenrm

\noindent{\sc Proof}: As in the proof of Corollary \top.11 this follows from 
Lemma \top.12 by cancelling the $t-1$ factors in $\grD(t)$ and evaluating at
$t=1.$ \hfill\za

\noindent{\sc Remark} \top.17: One can also write $d=m_4m_1$ or
$d=m_{2,1}m_{2,2}m_1$ where in each case the $m_s$ are pairwise relatively
prime positive integers. Also in the first case $d/w_i=m_4/v_i$ for 4 values of
$i$ and $d/w_j=m_1/v_j$ for the remaining index. In the second case
$d/w_i=m_{2,1}/v_i$ and $d/w_j=m_{2,2}/v_j$ for two pairs of index and
$d/w_k=m_1/v_k$ for the remaining index. In both cases one finds rational
homology spheres without any further conditions; however, one can also easily
show that the weights are not well-formed in either case.

\centerline{
\vbox{\tabskip=0pt \offinterlineskip
\def\tablerule{\noalign{\hrule}}
\halign to380pt {\strut#& \vrule#\tabskip=1em plus2em&
     \hfil#& \vrule#& \hfil#& \vrule#&  \hfil#& \vrule#&
     \hfil#& \vrule#\tabskip=0pt\cr\tablerule
\omit&height2pt&\multispan{7}&\cr
&&\multispan{7}\hfil {\bf Table}: 
\bbq-Homology 7-Spheres $M^7_{\bfw,d}$ admitting S-E
Structures\hfil&\cr\tablerule 
&&\omit\hidewidth $\bfw=(w_0,w_1,w_2,w_3,w_4)$\hidewidth&& 
\omit\hidewidth $d$\hidewidth&&
\omit\hidewidth $\mu$\hidewidth&&
\omit\hidewidth Order of $H_3(M^7_{\bfw,d},\bbz)$\hidewidth&\cr\tablerule
&&$(17,34,75,125,175)$&&$425$&&$4416$&&$582622237229761=17^{12}$&\cr\tablerule
&&$(17,238,381,635,889)$&&$2159$&&$16272$&&$118587876497=17^9$&\cr\tablerule
&&$(19,57,100,125,175)$&&$475$&&$3168$&&$16983563041=19^8$&\cr\tablerule
&&$(49,334,525,668,763)$&&$2338$&&$4995$&&$37259704=2^3\cdot
167^3$&\cr\tablerule
&&$(49,573,1862,2483,4393)$&&$9359$&&$36720$&&$282475249=7^{10}$&\cr\tablerule
&&$(50,65,73,73,105)$&&$365$&&$1152$&&$28398241=73^4$&\cr\tablerule
&&$(52,127,381,533,559)$&&$1651$&&$5040$&&$260144641=127^4$&\cr\tablerule
&&$(55,160,373,373,905)$&&$1865$&&$5932$&&$19356878641=373^4$&\cr\tablerule
&&$(87,558,687,1331,1331)$&&$3993$&&$4320$&&$1771561=11^6$&\cr\tablerule
&&$(93,459,780,1331,1331)$&&$3993$&&$4320$&&$1771561=11^6$&\cr\tablerule
&&$(93,1011,2298,3401,3401)$&&$10203$&&$13600$&&$11566801=19^2\cdot
179^2$&\cr\tablerule
&&$(97,1531,2201,2775,3253)$&&$9856$&&$9855$&&$9856=2^7\cdot 7\cdot
11$&\cr\tablerule
&&$(99,318,465,881,881)$&&$2343$&&$3520$&&$776161=881^2$&\cr\tablerule
&&$(101,439,559,579,619)$&&$2296$&&$2295$&&$2296=2^3\cdot 7\cdot
41$&\cr\tablerule
&&$(101,1597,1996,2695,3693)$&&$10081$&&$10080$&&$10081=17\cdot
593$&\cr\tablerule
&&$(101,1597,2096,2495,3793)$&&$10081$&&$10080$&&$10081=17\cdot
593$&\cr\tablerule
&&$(101,1697,2296,2695,4093)$&&$10881$&&$10880$&&$10881=3^3\cdot 13\cdot
31$&\cr\tablerule
&&$(103,1321,2337,2845,3251)$&&$9856$&&$9855$&&$9856=2^7\cdot 7\cdot
11$&\cr\tablerule
&&$(108,267,507,881,881)$&&$2643$&&$3520$&&$776161=881^2$&\cr\tablerule
&&$(109,1616,2047,2693,4417)$&&$10881$&&$10880$&&$10881=3^3\cdot 13\cdot
31$&\cr\tablerule
&&$(111,329,407,423,470)$&&$1739$&&$1728$&&$1369=37^2$&\cr\tablerule
&&$(111,658,2303,3071,6031)$&&$12173$&&$24600$&&$35611289=7^3\cdot
47^3$&\cr\tablerule
&&$(111,768,2523,3401,3401)$&&$10203$&&$13600$&&$11566801=19^2\cdot
179^2$&\cr\tablerule
&&$(113,1115,6021,8362,9589)$&&$25199$&&$50064$&&$1142897=113^3$&\cr\tablerule
&&$(113,1561,3345,8362,11819)$&&$25199$&&$50064$&&$1142897=113^3$&\cr\tablerule 
&&$(115,341,523,591,727)$&&$2296$&&$2295$&&$2296=2^3\cdot
7\cdot 41$&\cr\tablerule
&&$(115,797,949,987,2050)$&&$4897$&&$4896$&&$4897=59\cdot 83$&\cr\tablerule
&&$(125,1732,4577,5567,6433)$&&$18433$&&$18432$&&$18433$&\cr\tablerule
&&$(125,2599,4208,5569,9901)$&&$22401$&&$22400$&&$22401=3^2\cdot 19\cdot
131$&\cr\tablerule
&&$(127,1888,2643,4657,6671)$&&$15985$&&$15984$&&$15985=5\cdot 23\cdot
139$&\cr\tablerule
&&$(127,2266,3651,6043,8435)$&&$20521$&&$20520$&&$20521$&\cr\tablerule
&&$(127,2392,3399,6043,8561)$&&$20521$&&$20520$&&$20521$&\cr\tablerule
&&$(127,2770,4407,7429,10325)$&&$25057$&&$25056$&&$25057$&\cr\tablerule
&&$(129,511,1192,1235,1831)$&&$4897$&&$4896$&&$4897=59\cdot 83$&\cr\tablerule
&&$(133,346,379,527,857)$&&$2241$&&$2240$&&$2241=3^3\cdot 83$&\cr\tablerule
&&$(136,337,421,455,893)$&&$2241$&&$2240$&&$2241=3^3\cdot 83$&\cr\tablerule
&&$(136,2023,3237,5395,7553)$&&$18343$&&$17280$&&$289=17^2$&\cr\tablerule
&&$(137,1223,1427,2786,4349)$&&$9921$&&$9920$&&$9921=3\cdot 3307$&\cr\tablerule
&&$(137,1495,1699,3466,5301)$&&$12097$&&$12096$&&$12097$&\cr\tablerule
&&$(138,171,393,701,701)$&&$2103$&&$2800$&&$149401=701^2$&\cr\tablerule
&&$(139,2343,3721,6202,10061)$&&$22465$&&$22464$&&$22465= 5\cdot
4493$&\cr\tablerule}} } 

\vfill\eject
\centerline{
\vbox{\tabskip=0pt \offinterlineskip
\def\tablerule{\noalign{\hrule}}
\halign to380pt {\strut#& \vrule#\tabskip=1em plus2em&
     \hfil#& \vrule#& \hfil#& \vrule#&  \hfil#& \vrule#&
     \hfil#& \vrule#\tabskip=0pt\cr\tablerule
\omit&height2pt&\multispan{7}&\cr
&&\multispan{7}\hfil 
\bbq-Homology 7-Spheres $M^7_{\bfw,d}$ admitting S-E
Structures (cont.)\hfil&\cr\tablerule 
&&\omit\hidewidth $\bfw=(w_0,w_1,w_2,w_3,w_4)$\hidewidth&& 
\omit\hidewidth $d$\hidewidth&&
\omit\hidewidth $\mu$\hidewidth&&
\omit\hidewidth Order of $H_3(M^7_{\bfw,d},\bbz)$\hidewidth&\cr\tablerule
&&$(139,3171,5101,8548,13787)$&&$30745$&&$30744$&&$30745= 5\cdot
11\cdot 13\cdot 43$&\cr\tablerule
&&$(141,652,1909,2375,2701)$&&$7777$&&$7776$&&$7777= 7\cdot
11\cdot 101$&\cr\tablerule
&&$(141,1259,1492,3031,4663)$&&$10585$&&$10584$&&$10585= 5\cdot
29\cdot 73$&\cr\tablerule
&&$(141,2224,3475,4031,9729)$&&$19599$&&$19600$&&$19881=3^2\cdot
47^2$&\cr\tablerule
&&$(141,2363,3058,4309,9729)$&&$19599$&&$19600$&&$19881=3^2\cdot
47^2$&\cr\tablerule
&&$(143,194,209,231,291)$&&$1067$&&$1152$&&$9409=97^2$&\cr\tablerule
&&$(143,1135,2057,3476,5675)$&&$12485$&&$13608$&&$1288225=5^2\cdot
227^2$&\cr\tablerule
&&$(143,2981,6530,9795,16467)$&&$35915$&&$39168$&&$10660225=5^2\cdot
363^2$&\cr\tablerule
&&$(145,2157,3451,5752,9347)$&&$20851$&&$20160$&&$841=29^2$&\cr\tablerule
&&$(146,869,1955,2969,3983)$&&$9921$&&$9920$&&$9921=3\cdot 3307$&\cr\tablerule
&&$(147,207,230,245,299)$&&$1127$&&$1152$&&$2401=7^4$&\cr\tablerule
&&$(147,255,1056,1457,1457)$&&$4371$&&$5824$&&$2122849
=31^2\cdot 47^2$&\cr\tablerule
&&$(154,535,739,1427,2115)$&&$4969$&&$4968$&&$4969$&\cr\tablerule
&&$(154,763,1297,2975,3891)$&&$9079$&&$10368$&&$1682209=1297^2$&\cr\tablerule
&&$(155,921,1612,2149,4681)$&&$9517$&&$9792$&&$94249=307^2$&\cr\tablerule
&&$(155,1075,3532,5835,7064)$&&$17660$&&$21186$&&$12475024=3532^2$&\cr\tablerule
&&$(155,2309,7543,10006,12469)$&&$32481$&&$32480$&&$32481=3^4\cdot
401$&\cr\tablerule
&&$(155,2617,8467,11392,14163)$&&$36793$&&$36792$&&$36793$&\cr\tablerule
&&$(157,269,637,665,1090)$&&$2817$&&$2816$&&$2817=3^2\cdot 313$&\cr\tablerule
&&$(157,436,1401,1775,1993)$&&$5761$&&$5760$&&$5761=7\cdot 823$&\cr\tablerule
&&$(157,545,1051,1401,2608)$&&$5761$&&$5760$&&$5761=7\cdot 823$&\cr\tablerule
&&$(157,883,1558,2597,4311)$&&$9505$&&$9504$&&$9505=5\cdot
1901$&\cr\tablerule
&&$(157,1195,2182,3689,6027)$&&$13249$&&$13248$&&$13249$&\cr\tablerule
&&$(157,2339,5146,7641,12943)$&&$28225$&&$28224$&&$28225=5^2\cdot
1129$&\cr\tablerule
&&$(157,2651,5770,8733,14659)$&&$31969$&&$31968$&&$31969=7\cdot
4567$&\cr\tablerule
&&$(159,2365,3784,6307,12455)$&&$25069$&&$25488$&&$223729=11^2\cdot
43^2$&\cr\tablerule
&&$(163,1939,2747,6786,8887)$&&$20521$&&$20520$&&$20521$&\cr\tablerule
&&$(163,1939,3070,5171,10179)$&&$20521$&&$20520$&&$20521$&\cr\tablerule
&&$(163,2101,3394,5657,11151)$&&$22465$&&$22464$&&$22465=5\cdot
4493$&\cr\tablerule
&&$(166,1237,4371,5773,7175)$&&$18721$&&$18720$&&$18721=97\cdot
193$&\cr\tablerule
&&$(166,1399,2551,5513,7077)$&&$16705$&&$16704$&&$16705=5\cdot 13\cdot
257$&\cr\tablerule &&$(167,717,1324,1765,3805)$&&$7777$&&$7776$&&$7777=7\cdot
11\cdot 101$&\cr\tablerule
&&$(169,2015,6549,8732,10915)$&&$28379$&&$26208$&&$169=13^2$&\cr\tablerule
&&$(170,1267,3041,4477,7687)$&&$16641$&&$16640$&&$16641=3^2\cdot
43^2$&\cr\tablerule
&&$(171,247,556,695,973)$&&$2641$&&$2520$&&$361=19^2$&\cr\tablerule
&&$(175,271,299,306,751)$&&$1801$&&$1800$&&$1801$&\cr\tablerule
&&$(175,289,925,2312,3525)$&&$7225$&&$14688$&&$180124137569=17^6$&\cr\tablerule
&&$(175,2434,10605,13387,15995)$&&$42595$&&$43776$&&$1481089=1217^2$&\cr\tablerule
&&$(176,295,317,481,973)$&&$2241$&&$2240$&&$2241=3^3\cdot 83$&\cr\tablerule
&&$(176,1135,1397,4103,5675)$&&$12485$&&$13608$&&$1288225=5^2\cdot
227^2$&\cr\tablerule
&&$(177,997,3695,5044,6217)$&&$16129$&&$16128$&&$16129=3^6$&\cr\tablerule}} }
\vfill\eject \centerline{
\vbox{\tabskip=0pt \offinterlineskip
\def\tablerule{\noalign{\hrule}}
\halign to380pt {\strut#& \vrule#\tabskip=1em plus2em&
     \hfil#& \vrule#& \hfil#& \vrule#&  \hfil#& \vrule#&
     \hfil#& \vrule#\tabskip=0pt\cr\tablerule
\omit&height2pt&\multispan{7}&\cr
&&\multispan{7}\hfil 
\bbq-Homology 7-Spheres $M^7_{\bfw,d}$ admitting S-E
Structures (cont.)\hfil&\cr\tablerule 
&&\omit\hidewidth $\bfw=(w_0,w_1,w_2,w_3,w_4)$\hidewidth&& 
\omit\hidewidth $d$\hidewidth&&
\omit\hidewidth $\mu$\hidewidth&&
\omit\hidewidth Order of $H_3(M^7_{\bfw,d},\bbz)$\hidewidth&\cr\tablerule
&&$(177,2275,5950,7175,15399)$&&$30975$&&$30976$&&$31329=3^2\cdot
59^2$&\cr\tablerule
&&$(178,883,2031,3973,5033)$&&$12097$&&$12096$&&$12097$&\cr\tablerule
&&$(181,1613,4301,5197,11110)$&&$22401$&&$22400$&&$22401=3^2\cdot 19\cdot
131$&\cr\tablerule
&&$(181,1618,7011,8809,10607)$&&$28225$&&$28224$&&$28225=5^2\cdot
1129$&\cr\tablerule
&&$(181,2338,10251,12949,15467)$&&$41185$&&$41184$&&$41185=5\cdot 8237\cdot
131$&\cr\tablerule &&$(183,787,1453,2422,4661)$&&$9505$&&$9504$&&$9505=5\cdot
1901$&\cr\tablerule
&&$(183,1031,2608,4003,6793)$&&$14617$&&$14616$&&$14617=47\cdot
311$&\cr\tablerule 
&&$(185,477,991,1321,2788)$&&$5761$&&$5760$&&$5761=7\cdot
823$&\cr\tablerule
&&$(187,409,539,1320,2045)$&&$4499$&&$4896$&&$167281=409^2$&\cr\tablerule
&&$(187,631,1299,2746,3563)$&&$8425$&&$8424$&&$8425=5^2\cdot 337$&\cr\tablerule
&&$(187,781,2306,3459,5951)$&&$12683$&&$13824$&&$1329409=1153^2$&\cr\tablerule
&&$(187,1853,2594,6485,10931)$&&$22049$&&$23328$&&$1682209=1297^2$&\cr\tablerule
&&$(187,2416,8177,10965,19328)$&&$41072$&&$43470$&&$5837056=2^8\cdot
151^2$&\cr\tablerule 
&&$(191,235,433,509,904)$&&$2241$&&$2240$&&$2241=3^3\cdot
83$&\cr\tablerule
&&$(191,2467,7990,10627,18787)$&&$40041$&&$40040$&&$40041=3^3\cdot
1483$&\cr\tablerule
&&$(191,3607,11770,15757,27717)$&&$59041$&&$59040$&&$59041=17\cdot
23\cdot 151$&\cr\tablerule
&&$(193,3247,4202,10877,18335)$&&$36863$&&$36864$&&$37249=193^2$&\cr\tablerule
&&$(194,539,1155,2425,3157)$&&$7469$&&$7488$&&$9409=97^2$&\cr\tablerule
&&$(194,693,847,2425,3311)$&&$7469$&&$7488$&&$9409=97^2$&\cr\tablerule
&&$(195,484,613,1291,2387)$&&$4969$&&$4968$&&$4969$&\cr\tablerule
&&$(196,681,827,2383,3259)$&&$7345$&&$7344$&&$7345=5\cdot
13\cdot 113$&\cr\tablerule
&&$(196,1119,1411,4135,5741)$&&$12601$&&$12600$&&$12601$&\cr\tablerule
&&$(196,2337,7595,10127,17917)$&&$38171$&&$37440$&&$2401=7^4$&\cr\tablerule
&&$(197,881,4111,5188,6265)$&&$16641$&&$16640$&&$16641=3^2\cdot
43^2$&\cr\tablerule
&&$(197,1273,6071,7736,9205)$&&$24481$&&$24480$&&$24481$&\cr\tablerule
&&$(199,305,377,628,1309)$&&$2817$&&$2816$&&$2817=3^2\cdot 313$&\cr\tablerule
&&$(199,376,831,1603,2177)$&&$5185$&&$5184$&&$5185=5\cdot 17\cdot 
61$&\cr\tablerule
&&$(199,673,2811,3880,4751)$&&$12313$&&$12312$&&$12313=7\cdot
1759$&\cr\tablerule
&&$(199,1973,3157,7300,12429)$&&$25057$&&$25056$&&$25057$&\cr\tablerule
&&$(203,1409,1912,4931,8251)$&&$16705$&&$16704$&&$16705=5\cdot
13\cdot 257$&\cr\tablerule
&&$(205,389,1093,1389,2686)$&&$5761$&&$5760$&&$5761=7\cdot
823$&\cr\tablerule
&&$(205,1629,7127,10588,12421)$&&$31969$&&$31968$&&$31969=7\cdot
4567$&\cr\tablerule
&&$(206,247,259,319,771)$&&$1801$&&$1800$&&$1801$&\cr\tablerule
&&$(206,1331,4607,6143,10955)$&&$23241$&&$23240$&&$23241=3\cdot
61\cdot 127$&\cr\tablerule
&&$(207,2468,8021,10695,21183)$&&$42573$&&$43120$&&$380689=617^2$&\cr\tablerule
&&$(208,439,1007,2091,2737)$&&$6481$&&$6480$&&$6481$&\cr\tablerule
&&$(209,707,2038,3161,5407)$&&$11521$&&$11520$&&$11521=41\cdot
281$&\cr\tablerule
&&$(209,1351,5092,6859,12159)$&&$25669$&&$27000$&&$1825201=7^2\cdot
193^2$&\cr\tablerule  
&&$(211,339,436,985,1631)$&&$3601$&&$3600$&&$3601=13\cdot
277$&\cr\tablerule
&&$(211,1886,6077,8173,16135)$&&$32481$&&$32480$&&$32481=3^4\cdot
401$&\cr\tablerule
&&$(211,2306,7547,10063,19915)$&&$40041$&&$40040$&&$40041=3^3\cdot
1483$&\cr\tablerule}} }

\vfill\eject
\centerline{
\vbox{\tabskip=0pt \offinterlineskip
\def\tablerule{\noalign{\hrule}}
\halign to380pt {\strut#& \vrule#\tabskip=1em plus2em&
     \hfil#& \vrule#& \hfil#& \vrule#&  \hfil#& \vrule#&
     \hfil#& \vrule#\tabskip=0pt\cr\tablerule
\omit&height2pt&\multispan{7}&\cr
&&\multispan{7}\hfil  
\bbq-Homology 7-Spheres $M^7_{\bfw,d}$ admitting S-E
Structures (cont.)\hfil&\cr\tablerule 
&&\omit\hidewidth $\bfw=(w_0,w_1,w_2,w_3,w_4)$\hidewidth&& 
\omit\hidewidth $d$\hidewidth&&
\omit\hidewidth $\mu$\hidewidth&&
\omit\hidewidth Order of $H_3(M^7_{\bfw,d},\bbz)$\hidewidth&\cr\tablerule
&&$(214,489,659,1849,2551)$&&$5761$&&$5760$&&$5761=7\cdot
823$&\cr\tablerule
&&$(217,731,1075,2752,4557)$&&$9331$&&$9504$&&$47089=7^2\cdot
31^2$&\cr\tablerule
&&$(217,817,946,2795,4557)$&&$9331$&&$9504$&&$47089=7^2\cdot
31^2$&\cr\tablerule
&&$(217,2795,7310,13115,23219)$&&$46655$&&$46656$&&$47089=7^2\cdot
31^2$&\cr\tablerule
&&$(217,3440,5375,14405,23219)$&&$46655$&&$46656$&&$47089=7^2\cdot
31^2$&\cr\tablerule
&&$(220,237,1021,1367,1477)$&&$4321$&&$4320$&&$4321=29\cdot
149$&\cr\tablerule
&&$(221,317,805,1073,2194)$&&$4609$&&$4608$&&$4609=11\cdot 419$&\cr\tablerule
&&$(221,2416,5491,13617,19328)$&&$41072$&&$43470$&&$5837056=2^8\cdot
151^2$&\cr\tablerule
&&$(223,256,563,1041,1519)$&&$3601$&&$3600$&&$3601=13\cdot 277$&\cr\tablerule
&&$(223,2212,3539,9511,15261)$&&$30745$&&$30744$&&$30745=5\cdot
11\cdot 13\cdot 43$&\cr\tablerule
&&$(223,2437,13071,18166,20825)$&&$54721$&&$54720$&&$54721$&\cr\tablerule
&&$(223,4211,9087,22606,31915)$&&$68041$&&$68040$&&$68041$&\cr\tablerule
&&$(226,2245,12123,16837,19307)$&&$50737$&&$50400$&&$12769=113^2$&\cr\tablerule
&&$(226,3143,6735,16837,23797)$&&$50737$&&$50400$&&$12769=113^2$&\cr\tablerule
&&$(227,721,856,2523,4099)$&&$8425$&&$8424$&&$8425=5^2\cdot
337$&\cr\tablerule
&&$(227,901,1051,3228,5179)$&&$10585$&&$10584$&&$10585=5\cdot
29\cdot 73$&\cr\tablerule
&&$(227,1015,3496,4737,9247)$&&$18721$&&$18720$&&$18721=97\cdot
193$&\cr\tablerule
&&$(227,1241,4400,5867,11507)$&&$23241$&&$23240$&&$23241=3\cdot
61\cdot 127$&\cr\tablerule
&&$(229,2503,5461,13652,19341)$&&$41185$&&$41184$&&$41185=5\cdot
8237$&\cr\tablerule
&&$(235,323,334,891,1459)$&&$3241$&&$3240$&&$3241=7\cdot 463$&\cr\tablerule
&&$(237,275,766,1021,2023)$&&$4321$&&$4320$&&$4321=29\cdot
149$&\cr\tablerule
&&$(238,1301,7451,10289,11827)$&&$31105$&&$31104$&&$31105=5\cdot
6221$&\cr\tablerule
&&$(241,639,3275,4792,5671)$&&$14617$&&$14616$&&$14617=47\cdot
311$&\cr\tablerule
&&$(242,385,409,1419,2045)$&&$4499$&&$4896$&&$167281=409^2$&\cr\tablerule
&&$(243,457,2339,3280,3979)$&&$10297$&&$10296$&&$10297=7\cdot
1471$&\cr\tablerule
&&$(245,2434,5355,13387,21175)$&&$42595$&&$46776$&&$1481089=1217^2$&\cr\tablerule
&&$(247,292,799,1583,2121)$&&$5041$&&$5040$&&$5041=71^2$&\cr\tablerule
&&$(247,1351,3439,8474,12159)$&&$25669$&&$27000$&&$1825201=1361^2$&\cr\tablerule
&&$(247,3190,6871,17177,27237)$&&$54721$&&$54720$&&$54721$&\cr\tablerule
&&$(250,275,393,917,1441)$&&$3275$&&$3168$&&$625=5^4$&\cr\tablerule
&&$(250,1367,3231,8077,11557)$&&$24481$&&$24480$&&$24481$&\cr\tablerule
&&$(251,561,592,1963,3115)$&&$6481$&&$6480$&&$6481$&\cr\tablerule
&&$(253,311,1096,1377,2725)$&&$5761$&&$5760$&&$5761=7\cdot
823$&\cr\tablerule
&&$(253,1507,6530,9795,17831)$&&$35915$&&$39168$&&$10660225=5^2\cdot
653^2$&\cr\tablerule
&&$(255,844,1519,4135,6497)$&&$13249$&&$13248$&&$13249$&\cr\tablerule
&&$(259,285,407,950,1615)$&&$3515$&&$3456$&&$1369=37^2$&\cr\tablerule
&&$(259,643,1993,3536,6171)$&&$12601$&&$12600$&&$12601$&\cr\tablerule
&&$(261,287,410,957,1653)$&&$3567$&&$3520$&&$1681=41^2$&\cr\tablerule
&&$(261,491,1762,2773,4795)$&&$10081$&&$10080$&&$10081=17\cdot
593$&\cr\tablerule &&$(262,291,331,883,1475)$&&$3241$&&$3240$&&$3241=7\cdot
463$&\cr\tablerule &&$(262,443,469,1641,2371)$&&$5185$&&$5184$&&$5185=5\cdot
17\cdot 61$&\cr\tablerule}} }
\vfill\eject
\centerline{
\vbox{\tabskip=0pt \offinterlineskip
\def\tablerule{\noalign{\hrule}}
\halign to380pt {\strut#& \vrule#\tabskip=1em plus2em&
     \hfil#& \vrule#& \hfil#& \vrule#&  \hfil#& \vrule#&
     \hfil#& \vrule#\tabskip=0pt\cr\tablerule
\omit&height2pt&\multispan{7}&\cr
&&\multispan{7}\hfil 
\bbq-Homology 7-Spheres $M^7_{\bfw,d}$ admitting S-E
Structures (cont.)\hfil&\cr\tablerule 
&&\omit\hidewidth $\bfw=(w_0,w_1,w_2,w_3,w_4)$\hidewidth&& 
\omit\hidewidth $d$\hidewidth&&
\omit\hidewidth $\mu$\hidewidth&&
\omit\hidewidth Order of $H_3(M^7_{\bfw,d},\bbz)$\hidewidth&\cr\tablerule
&&$(263,1699,3921,9802,15421)$&&$31105$&&$31104$&&$31105=5\cdot
6221$&\cr\tablerule
&&$(268,511,911,2599,3777)$&&$8065$&&$8064$&&$8065=5\cdot
1613$&\cr\tablerule
&&$(271,377,673,1696,2745)$&&$5761$&&$5760$&&$5761=7\cdot
823$&\cr\tablerule
&&$(275,729,2734,4465,7927)$&&$16129$&&$16128$&&$16129=127^2$&\cr\tablerule
&&$(277,441,2591,3748,4465)$&&$11521$&&$11520$&&$11521=41\cdot
281$&\cr\tablerule
&&$(299,325,1869,2492,3115)$&&$8099$&&$7488$&&$169=13^2$&\cr\tablerule
&&$(301,363,1294,2257,3851)$&&$8065$&&$8064$&&$8065=5\cdot
1613$&\cr\tablerule
&&$(301,3289,17342,24219,44849)$&&$89999$&&$90000$&&$90601=7^2\cdot
43^2$&\cr\tablerule
&&$(301,3887,12259,28704,44849)$&&$89999$&&$90000$&&$90601=7^2\cdot
43^2$&\cr\tablerule
&&$(307,2441,13120,18307,33867)$&&$68041$&&$68040$&&$68041$&\cr\tablerule
&&$(309,349,2279,3244,3901)$&&$10081$&&$10080$&&$10081=17\cdot
593$&\cr\tablerule
&&$(311,495,2104,3403,6001)$&&$12313$&&$12312$&&$12313=7\cdot
1759$&\cr\tablerule
&&$(311,2473,8037,18856,29365)$&&$59041$&&$59040$&&$59041=17\cdot 23\cdot
151$&\cr\tablerule
&&$(316,1727,9577,13345,24648)$&&$49612$&&$49770$&&$49856=2^4\cdot
79^2$&\cr\tablerule
&&$(316,2041,6751,15857,24648)$&&$49612$&&$49770$&&$49856=2^4\cdot
79^2$&\cr\tablerule
&&$(328,347,449,1571,2347)$&&$5041$&&$5040$&&$5041=71^2$&\cr\tablerule
&&$(339,383,1780,2839,4957)$&&$10297$&&$10296$&&$10297=7\cdot
1471$&\cr\tablerule
&&$(341,407,2306,3459,6171)$&&$12683$&&$13824$&&$1329409=1153^2$&\cr\tablerule
&&$(356,387,2225,2967,5547)$&&$11481$&&$11440$&&$7921=89^2$&\cr\tablerule
&&$(357,388,2231,2975,5593)$&&$11543$&&$11520$&&$9409=97^2$&\cr\tablerule}} }

\bigskip
\bigskip
\centerline{\bf \dis. Discussion of the Table}   
\bigskip

In this section we give a discussion about the representatives listed in the
Table. It is easy to notice the existence of {\it twins} in the Table. These 
are rational homology 7-spheres with the same degree $d$, Milnor number $\mu$
and order of $H_3.$ Twins often occur as adjacent listings with the same
$w_0,$ but this is not always the case as with twins $d=|H_3|=10881,\mu=10880$ 
and $w_0=101$ and $109,$ and the twins $d=|H_3|=7777$ with $w_0=141$ and 
$w_0=167.$ Twins may also be members of a larger set, such as the {\it 
septuplets} with $d=|H_3|=5761$ and $\mu=5760.$ These have 
$w_0=157,157,185,205,214,253,271,$ respectively. Since twins have the same 
Milnor number, it is tempting to conjecture that twins correspond to 
homeomorphic or even diffeomorphic links, but we have no proof as of yet. In 
fact, except for cases where the order of $H_3$ contains no primes of order 
larger than one in its prime decomposition, we don't even know that twins have 
isomorphic $H_3$'s. Notice that the order of $H_3$ tends to be quite large 
varying from $169=13^2$ to $17^{12}$ a number over $500$ trillion. 

Another interesting fact is that of the 184 rational homology 7-spheres listed 
in the Table, only 10 have even degree, while the remaining 174 have odd 
degree, and the degree is even if and only if the order of $H_3$ is even. But 
even more intriguing is the fact that for all 174 rational homology 7-spheres 
with odd degree, the order $|H_3|\equiv1~\hbox{mod}(8).$ In [BGN4] we 
construct positive Sasakian structures on homotopy 9-spheres using the 
rational homology 7-spheres listed in the Table. There we show that the exotic 
Kervaire sphere can only occur when the degree of the rational homology sphere 
is even. 

Also of interest are invariants of the underlying contact, and almost contact 
structures. The underlying almost contact structures are classified [Sa] by 
homotopy classes of maps $[M^7,SO(8)/U(4)],$ and Morita [Mo] shows that for 
Brieskorn spheres this is a function of the Milnor number $\mu.$ It seems 
reasonable that a similar result holds true in our case. There are candidates 
for this in the table. For example the rational homology 7-spheres with 
weights $\bfw=(196,2337,7595,10127,17917)$ and degree $d=38171,$ and with 
weights $\bfw=(147,207,230,245,299)$ and degree $d=1127$ both have 
$|H_3|=7^4,$ so they could be diffeomorphic. But they have very different 
Milnor numbers, namely, $37440$ and $1152,$ respectively, so they could belong 
to distinct almost contact structures. Similarly there are 4 rational homology 
7-spheres with $|H_3|=97^2,$ two are twins having the same Milnor number, but 
the other two have different Milnor numbers both different than the Milnor 
number of the twins. Moreover, twins probably belong to the same underlying 
almost contact structures, but could possibly belong to distinct contact 
structures. It appears that nothing is known beyond homotopy spheres [Us1, 
Us2] about distinct contact structures within the same underlying almost 
contact structures. 

\bigskip
\bigskip
\centerline{\bf \reg. Some Comments on Regular Rational Homology Spheres}   
\bigskip

In this section we discuss some rational homology spheres that are regular, in
particular the homogeneous ones. The following result follows easily from
previous work [BG1] together with the well-known classification of del Pezzo
surfaces:

\noindent{\sc Proposition} \reg.1: \tensl Let $\cals=(g,\xi,\eta,\Phi)$ be a
regular positive Sasakian structure on a smooth compact 5-manifold $M^5.$
Then $M^5$ is a rational homology sphere if and only if $M^5$ is covered by
$S^5$ and $\cals$ is homologous to the standard Sasakian structure with
the round metric $g_0$. \tenrm

It is well-known that the standard Sasakian structure is a homogeneous \Se
structure.  Dimension seven is a bit more interesting:

\noindent{\sc Theorem} \reg.2: \tensl Let $\cals=(g,\xi,\eta,\Phi)$ be a
regular positive Sasakian structure on a smooth compact 7-manifold $M^7.$ Then
$M^7$ is a rational homology sphere if and only if it is one of the following:
\item{1.} $M^7=S^7$ and $\cals$ is homologous to the standard Sasakian
structure with the round metric.
\item{2.} $M^7=V_2(\bbr^5)$ the Stiefel manifold of $2$-frames in $\bbr^5$
and $\cals$ is homologous to the standard homogeneous \Se structure on
$V_2(\bbr^5)$ [BG1,BG2].
\item{3.} $M^7$ is a circle bundle over a smooth variety $V_5$ of degree $5$ in
$\bbc\bbp^6$ with a compatible Sasakian structure $\cals.$
\item{4.} $M^7$ is a circle bundle over a smooth variety $V_{22}$ of degree
$22$ in $\bbc\bbp^{13}$ with a compatible Sasakian structure $\cals.$

\noindent Furthermore, $M^7$ admits a homogeneous \Se structure if and only if
$M^7=S^7$ or $V_2(\bbr^5).$  \tenrm

\noindent{\sc Proof}: By [BG1] $M^7$ is a regular rational homology sphere
with a Sasakian structure $\cals$ if and only if it is the total space of an
$S^1$ bundle over a smooth projective 3-fold $\calz$ with the same rational
homology groups as projective space $\bbc\bbp^3.$ Furthermore,  $\cals$ is
positive [BGN3] if and only if $\calz$ is Fano. Thus, $\calz$ must occur on
Iskovskikh's list [Isk] (see Remark \reg.3 below) of smooth Fano 3-folds of 
the first kind, and there are precisely four which have the same rational 
cohomology groups as $\bbc\bbp^3.$ This gives the four cases above. The last 
statement follows from Corollary 4.1.3 of [BG2]. \hfill\za 

\noindent{\sc Remarks} \reg.3: (1) Case 4 in Theorem \reg.2 has an interesting 
history. The 3-fold $V_{22}$ was missed by Fano in his original classification 
of smooth 3-folds with an ample anti-canonical line bundle. It was then found 
by Iskovskikh [Ish] in his study of Fano's work, but a mistake was made and 
not all were found. Mukai and Umemura [MU] (See also [IsPr]) produced a 
$V_{22}$ that is an equivariant compactification of $SL(2,\bbc)/\bbi$ that was 
missed by Iskovskikh.  Here $\bbi$ is the icosahedral group. Later Prokhorov 
(see Proposition 4.3.11 of [IsPr]) showed that the Mukai-Umemura $V_{22}$ 
completes the Fano-Iskovskikh classification of Fano 3-folds. Recently Tian 
[Ti1,Ti2] showed that there are  deformations $P_a$ of the Mukai-Umemura 
$V_{22}$ which do not admit a K\"ahler-Einstein structure, giving a 
counterexample to the folklore conjecture that every every compact K\"ahler 
manifold with no holomorphic vector fields admits a compatible 
K\"ahler-Einstein metric. Thus, 
the Sasakian circle bundle over $P_a$ does not admit a compatible \Se metric.
(2) In the four cases of Theorem \reg.2, the corresponding Fano 3-folds are 
precisely those Fano 3-folds that are almost homogeneous with respect to the 
group $SL(2,\bbc).$ (See [IsPr], pg 116). 

There is a straightforward procedure for finding all rational homology spheres
$M^{2n+1}$ that admit a homogeneous \Se structure. By Theorem 3.2 of [BG1]
$M^{2n+1}$ must fibre over a generalized flag manifold $G/P,$ where $G$ is a
complex semi-simple Lie group, and $P$ is a parabolic subgroup. In order that
$M^{2n+1}$ be a rational homology sphere, it is necessary that $G/P$ have the
rational homology of a projective space. Hence, we may restrict ourselves to
the case where $G$ is simple and $P$ is maximal parabolic. The procedure for
computing the cohomology ring of $G/P$ is outlined in Baston and Eastwood
[BE]. All $G/P$'s with $G$ simple are realized by crossing out nodes in each
Dynkin diagram of $G.$ When $P$ is maximal parabolic only one node is crossed
out. The rank of the cohomology groups is determined by the Hasse diagram
$W^\gp$ which is the coset space $W_\gg/W_\gp$ where $W_\gg$ is the Weyl group
of the Lie algebra $\gg$ of $G,$ and $W_\gp$ is the Weyl group of the Levi
factor of the Lie algebra $\gp$ of $P.$ Then the cohomology groups of $G/P$
will have the same rank as $\bbc\bbp^n$ if and only if $W^\gp$ has precisely
one element of length $l$ for each $l=1,\cdots,n.$ One then needs to check all
maximal parabolics for all Dynkin diagrams, and compute the Hasse diagram for
each case. There are many cases and repetitions can and do occur. Here we 
mention the Stiefel manifolds $V_2(\bbr^{2n+1})$ which are circle bundles over
the odd quadrics $Q_{2n-1}$ and the homogeneous 3-Sasakian rational homology
sphere $G_2/Sp(1)_+$ (cf. [BGP] and Remark \reg.6(2) below). A Gysin sequence 
or spectral sequence argument shows that
$$H^p(V_2(\bbr^{2n+1}),\bbz)\approx \cases{\bbz &if $p=0,4n-1$;\cr
                                                          \bbz_2&if $p=2n$;\cr 
                                        0 &otherwise.\cr}$$ 
$$H^p(G_2/Sp(1)_+,\bbz)\approx \cases{\bbz &if $p=0,11$;\cr
                                                          \bbz_3&if $p=4,8$;\cr
                                         0 &otherwise.\cr}\leqno{\reg.4}$$
We have

\noindent{\sc Proposition} \reg.5: \tensl The Stiefel manifold
$V_2(\bbr^{2n+1})$ and $G_2/Sp(1)_+$ are simply connected rational homology
spheres which admit homogeneous \Se structures. 
\tenrm

\noindent{\sc Remarks} \reg.6: (1) Since $V_2(\bbr^{2n+1})$ can be represented 
as the link of the quadric hypersurface singularity, \reg.3 can be derived from
the Milnor-Orlik algorithm described in section 3. (2) There are two 
non-conjugate $Sp(1)$ subgroups of the exceptional Lie group $G_2,$ denoted in 
[BGP] as $Sp(1)_\pm.$ The quotient $G_2/Sp(1)_+$ has a homogeneous 3-Sasakian 
structure, whereas the quotient $G_2/Sp(1)_-$ does not. It does, however, have 
a homogeneous \Se structure, and as homogeneous \Se manifolds 
$G_2/Sp(1)_-\approx V_2(\bbr^7).$

There is an obvious corollary of Theorem 4.2.6 and Proposition 5.4.4 of [BG2], 
viz.

\noindent{\sc Corollary} \reg.7: \tensl Let $M^{4n+3}$ be a rational homology
sphere that admits a 3-Sasakian homogeneous structure. Then $M^{4n+3}$ is
either $S^{4n+3},\bbr\bbp^{4n+3},$ or $G_2/Sp(1)_+.$ \tenrm 


\bigskip
\centerline{\bf Appendix}
\bigskip
In this appendix, we generalize a result of Johnson and Koll\'ar (Proposition 
11 and Corollary 13 of [JK1]) to arbitrary orbifold Fano index $I.$ While this 
generalization is straightforward, we give more detail than [JK1].
We refer to our previous work [BGN1] and the book [KM] for background material.

Suppose $X \subset \bbp(\bfw) = \bbp(w_0,\ldots,w_n)$ 
is an irreducible hypersurface in 
weighted projective space.  We assume that $X$ has at worst quotient 
singularities (at singular points of $\bbp(\bfw)$) and that $X$ is Fano.  We 
would
like to study under what conditions $(X,D)$ is klt whenever $D \equiv
-\alpha K_X$ where $\alpha > {n-1\over n}$ and $D$ is effective.  
Under these circumstances $X$ admits a K\"ahler-Einstein metric of positive 
scalar curvature [DK,JK1]. We shall prove

\noindent{\sc Theorem} A.1: \tensl
Suppose $X \subset \bbp(\bfw) = \bbp(w_0,\ldots,w_n)$ is a normal Fano variety 
 of index $I$ satisfying
$$
\deg(X) < {nw_0w_1\over (n-1)I}.
$$
Then there exists $\epsilon > 0$ such that
$(X,D)$ is klt for every effective $D \equiv -{n-1+\epsilon \over n}K_X$.
\tenrm

\noindent{\sc Proof}:  By induction on dimension, which reduces to the surface
case already handled in [JK1] and [BGN1].  Suppose for a moment that $X$ is
a surface with at worst quotient singularities.  
Johnson and Kollar [JK1] give sufficient conditions for $(X,D)$ to be
klt:
\item{1.} If $D = \sum_{i=1}^r \alpha_iC_i$ then $\alpha_i < 1$ for all i,
\item{2.} for all smooth points $x \in X$, $\mult_x(D) \leq 1$,
\item{3.} if $P \subset X$ is a singular point and $\pi: Y \rightarrow X$
a local finite cover resolving the singularity at $P$
 then $\pi^\ast D$ has multiplicity at most one at $Q = \pi^{-1}(P)$.

Condition 1 is also necessary, while neither 2 nor 3 is necessary.  
We will now formulate analogues of 1,2, and 3 above, for an arbitrary
$(X,D)$, 
designed to guarantee that 1,2, and 3 will hold when we cut $X$ with the
appropriate number of hypersurfaces, allowing us to apply the inversion of
adjunction formula to conclude that $(X,D)$ is klt.

\item{(i)} If $D = \sum_{i=1}^r \alpha_iD_i$ then $\alpha_i < 1$ for all i, and
\item{(ii)} for all smooth points $x \in X$, $\mult_x(D) \leq 1$, and
\item{(iii)} suppose $x \in X$ is a singular point of $X$ with local group
$G_x$ and that
$H_1, H_2, \ldots, H_{n-1}$ are general hypersurfaces
through $x$; then we ask that the intersection number at $x,$
$ i(x,X \cdot H_1 \cdot \ldots \cdot H_{n-1};\bbp(\bfw))$ satisfies
$$ i(x,X \cdot H_1 \cdot \ldots \cdot H_{n-1};\bbp(\bfw)) \leq {1\over |G_x|}.
$$

Suppose that conditions i, ii, and iii are satisfied by a divisor $D$ on $X$
and suppose $x \in X$.  Choose general
hypersurfaces $H_1,H_2,\ldots, H_{n-1}$ through $x$ and 
let $V = X \cap H_1 \cap \ldots \cap H_{n -3}$ with $D_V = D \cap V$.
Consider the pair
$(V,D_V)$.  We see that $V$ is a surface and we
will show that $(V,D_V)$ satisfies 1,2, and 3 above and hence
$(V,D_V)$ is klt at $x$.  
Using the inversion of adjunction, $(X,D)$ is klt at $x$ as well.
  
Write $D_V = \sum b_jC_j$. 
Since
the hypersurfaces $H_i$ are general, we can use  Remark 8.2 of  [Ful] to
see that $b_j < 1$ for all $j$ as each $b_j$
is equal to one of the $\gra_i$.  Similarly, the multiplicity of $D$ at $x$ 
will be preserved under intersection by general hyperplanes and so
$\mult_x(D_V) \leq 1$ (see [Ful] Corollary 12.4).  
Finally, under hypothesis iii above, if $\pi: Y \rightarrow V$ is a local
cover of the quotient singularity at $x$, we see, using [Ful] 8.3.12, that
3 is satisfied.  More precisely, letting $y = \pi^{-1}(x)$,

$$\eqalign{
\mult_y(\pi^\ast D_V) &\leq \pi^\ast D_V \cdot \pi^\ast H_{n -2} \cdot
\pi^{\ast} H_{n-1} \cr
& \leq  |G_x|i(x,
D_V \cdot  H_{n -2} \cdot  H_{n-1};\bbp(\bfw)) \cr
&\leq   1,\cr}$$
the last inequality coming from iii.  
Thus it is sufficient, in order to prove Theorem 0.1,
 to verify conditions i, ii, and iii above where $D \equiv 
-{n-1+\epsilon\over n}K_X$ is an effective divisor.

Suppose then that we write
$$
D = \sum_{i=1}^r \alpha_iD_i.
$$
We will show that $\alpha_1 < 1$, the other cases being identical.  
Suppose $x \in D_1$ is a smooth point of $D_1$ and $X$; here we use the
fact that $X$ is normal and hence smooth along $D_1$.
To simplify notation in what follows, we 
let $E = \alpha_1 D_1$.  
We will assume for simplicity, rearranging the coordinates if necessary,
that 
$$
x = (x_0,\ldots,x_k,0,\ldots,0).
$$
Consider
the hyperplane $H_1$ given by $z_n = 0$.  If $E \subset H_1$ then replace
$H_1$ with $H_2$, given by $z_{n-1} = 0$ and start over.  
Assuming that $H_1$ meets $E$ properly,
write
$$
H_1 \cdot E = \sum \beta_jV_j.
$$
We then repeat the above procedure, replacing $E$ with $H_1 \cdot E$ and 
intersecting with $H_2$.  
We continue up to and including $H_{n-k}$ given by 
$z_{k+1} = 0$.
Next consider the hypersurface $D_i$ for $1 \leq i \leq k$ defined by 
$$
c_{i-1}z_{i-1}^{w_{i}} - c_{i}z_{i}^{w_{i-1}}
$$
where the $c_i$ are chosen so that $x \in D_i$.  
The divisors $H_1,\ldots,H_{n-k},D_1,\ldots,D_{k}$ cut out the point $x$ set 
theoretically. Thus we can continue the intersection process above, using
the $\bbq$-divisors $D_i/{\rm min}(w_{i-1},w_i)$, until
we  obtain a cycle $Z$ of dimension zero.  Bounding the degree of $Z$ by
the degree of the total intersection class we have
$$\eqalign{
\deg(Z) &\leq w_n\ldots w_3 \deg(E) \cr
        & \leq  w_n \ldots w_3 \deg(D),\cr}$$
since the worst case scenario for the degree of $Z$ occurs
when all possible intersections are proper and where no reordering of
variables has been necessary; note that the degrees in this formula are
relative to $\calo_{\bbp(\bfw)}(1)$.  
On the other hand since $\mult_x(H_i) \geq 1$
for all $i$ and  
$${\mult_x(D_i)\over {\rm min}(w_{i-1},w_i)} \geq 1,
$$
it follows by [Ful] Corollary 12.4 that
$$
\deg(Z) \geq  \mult_x(E).
$$
Hence
$$
\mult_x(E) \leq w_n\ldots w_3 \deg(D).
$$
But 
$$
\deg(D)  = {n-1+\epsilon\over n \prod_{i=0}^n w_i } I \deg(X),
$$
where $\deg(X)$ denotes the homogeneous degree of the polynomial defining $X$. 
So if we choose $\epsilon$ so that 
$$
{n-1+\epsilon\over n}I\deg(X)  < w_0w_1
$$
we see that 
$$
\alpha_1 \mult_x(D_1) =  \mult_x(E) < 1/w_2 \leq 1.
$$
Since we have assumed that
$X$  and $D_1$ are smooth at $x$ it follows that 
$\mult_x(D_1) = 1$ and hence $\alpha_1 < 1$ as desired.

\medskip

Next we address property ii.   The argument is essentially
identical to the verification of property i.  In particular, let 
$$
x = (x_0,x_1,\ldots,x_n);
$$
again we will assume that $x_0,\ldots, x_k$ are non--zero and that the
rest of the coordinates are zero, possibly after reordering the coordinates.  
Intersecting with the divisors $H_i$ and ${D_j\over w_j}$ as above will
show that
$$
\mult_x(D) < {1\over w_2} \leq 1.
$$
Again, this result holds with the original $w_2$ since the initial ordering
of the $w_i$'s is increasing.

Finally we establish property iii.  We begin with a specific example so that
the coordinates are simple and the argument transparent.  Suppose
$x = (0,\ldots,0,1) \in X$.  Thus $|G_x| \leq w_n$.  To see that
$$
i(x,X \cdot H_1 \cdot \ldots \cdot H_{n-1};\bbp(\bfw)) \leq {1\over w_n}
$$
we consider the intersection theoretic argument above.  In this case, $x$
is cut out by $H_i = \{z_i = 0\}$ for $0 \leq i \leq n-1$.  Thus we
never need to intersect with $\{z_n = 0\}$ and consequently instead of 
$$
\deg(Z) \leq w_n\ldots w_3 \deg(D),
$$
we will find
$$
\deg(Z) \leq w_{n-1}\ldots w_2 \deg(D).
$$
Writing $Z = m[x] + Z^\prime$ where $Z^\prime$ is a zero cycle not supported
at $x$, we find, plugging in the value of $\deg(D)$, that 
$$
m \leq {1\over w_n}.
$$
But $i(x,X\cdot H_1 \cdot \ldots \cdot H_{n-1};\bbp(\bfw))$ is the minimal 
intersection
number supported on $x$ which one can obtain but cutting out $x$ and thus
$$
i(x,X\cdot H_1 \cdot \ldots \cdot H_{n-1};\bbp(\bfw)) \leq {1\over w_n}
$$
as desired.

More generally, if
$P $ is some other singular point of $X$, with coordinates $(x_0,\cdots,x_n)$ 
then some of the homogeneous coordinates $\{z_i\}$
must be zero.  If
$j$ is the smallest index such that $z_j \neq 0$ then $|G_x| \leq w_j$.  
The argument then procedes as in case ii, the $w_j$ in the denominator
coming from taking the worst case scenario in computing the intersection
product. \hfill\za

Note that the most difficult of the three conditions to satisfy is definitely
the third.  This is the same as in the surface case but the higher dimensional
case is more difficult to deal with as, in the surface case, condition (iii) 
was dealt with using the inversion of the adjunction formula.  Again as in the 
surface case, we can weaken the condition somewhat if $X$ does not contain 
certain planes.

\medskip

\noindent{\sc Acknowledgments}: We would like to thank Jennifer Johnson and 
J\'anos Koll\'ar for making their computer list available, as well as for 
valuable comments and their interest in our work. We would also like to thank 
Alexi Kobalev and Gang Tian for helpful discussions. The second
author would like to thank Max-Planck-Institute f\"ur Mathematik in Bonn for
support and hospitality during the summer of 2001 when this paper was being 
completed.

\bigskip
\medskip
\centerline{\bf Bibliography}
\medskip
\medskip
\medskip
\medskip
\font\ninesl=cmsl9
\font\bsc=cmcsc10 at 10truept
\parskip=1.5truept
\baselineskip=11truept
\ninerm
\item{[BE]} {\bsc R.J. Baston and M.G. Eastwood}, {\ninesl The Penrose
Transform}, Oxford University Press, New York, 1989.
\item{[Be]} {\bsc A. Besse}, {\ninesl Einstein manifolds},
Springer-Verlag, Berlin-New York, 1987.
\item{[BFGK]} {\bsc H. Baum, T. Friedrich, R. Grunewald, and I. Kath},  
{\ninesl Twistors and Killing Spinors on Riemannian Manifolds}, 
Teubner-Texte f\"ur Mathematik, vol. 124, Teubner, Stuttgart, Leipzig, 1991. 
\item{[BG1]} {\bsc C. P. Boyer and  K. Galicki}, {\ninesl On Sasakian-Einstein
Geometry}, Int. J. Math. 11 (2000), 873-909.
\item{[BG2]} {\bsc C. P. Boyer and  K. Galicki}, {\ninesl 
3-Sasakian manifolds}. {\it Surveys in differential geometry: 
essays on Einstein manifolds}, 123--184, Surv. Differ. Geom.,
VI, Int. Press, Boston, MA, 1999.
\item{[BG3]} {\bsc C. P. Boyer and  K. Galicki}, {\ninesl New Einstein Metrics
in Dimension Five},  submitted for publication; math.DG/0003174.
\item{[BGN1]} {\bsc C. P. Boyer, K. Galicki, and M. Nakamaye}, {\ninesl On the
Geometry of Sasakian-Einstein 5-Manifolds}, submitted
for publication; math.DG/0012041.
\item{[BGN2]} {\bsc C. P. Boyer, K. Galicki, and M. Nakamaye}, {\ninesl On
Positive Sasakian Geometry}, submitted
for publication; math.DG/0104126.
\item{[BGN3]} {\bsc C. P. Boyer, K. Galicki, and M. Nakamaye}, {\ninesl
Sasakian-Einstein Structures on $\scriptstyle{9\#(S^2\times S^3)}$}, submitted
for publication; math.DG/0102181.
\item{[BGN4]} {\bsc C. P. Boyer, K. Galicki, and M. Nakamaye}, {\ninesl On the 
Sasakian Geometry of Homotopy Spheres}, to appear
\item{[BGP]} {\bsc C. P. Boyer, K. Galicki, and P. Piccinni}, {\ninesl 
3-Sasakian Geometry, Nilpotent Orbits, and Exceptional Quotients}, preprint 
DG/0007184 , to appear in Ann. Global Anal. Geom.
\item{[DK]} {\bsc J.-P. Demailly and J. Koll\'ar}, {\ninesl Semi-continuity of
complex singularity exponents and K\"ahler-Einstein metrics on Fano
orbifolds}, preprint AG/9910118, to appear in Ann. Scient. Ec. Norm. Sup. Paris
\item{[Dol]} {\bsc I. Dolgachev}, {\ninesl Weighted projective varieties}, in
Proceedings, Group Actions and Vector Fields, Vancouver (1981) LNM 956, 34-71.
\item{[Fle]} {\bsc A.R. Fletcher}, {\ninesl Working with weighted complete
intersections}, Preprint MPI/89-95, revised version in  {\it Explicit
birational geometry of 3-folds},  A. Corti and M. Reid, eds.,
Cambridge Univ. Press, 2000,  pp 101-173.
\item{[Ful]} {\bsc W. Fulton}, {\ninesl Intersection Theory}, Springer-Verlag, 
New York, 1984.
\item{[GS]} {\bsc K. Galicki and S. Salamon}, {\ninesl On
Betti numbers of 3-Sasakian manifolds},  Geom. Ded. 63 (1996), 45-68.
\item{[HZ]} {\bsc F. Hirzebruch and D. Zagier}, {\ninesl The Atiyah-Singer
Theorem and Elementary Number Theory}, Publish or Perish, Inc., Berkeley, 1974.
\item{[Isk]} {\bsc V.A. Iskovskikh}, {\ninesl Anticanonical Models of
Three-dimensional Algebraic Varieties}, J. Soviet Math. 13 (1980) 745-814.
\item{[IsPr]} {\bsc V.A. Iskovskikh and Yu.G. Prokhorov}, {\ninesl Fano 
Varieties}, Enc. Math. Sci. Vol 47, Algebraic Geometry V, A.N. Parshin and 
I.R. Shaferevich, Eds., Springer-Verlag, 1999.
\item{[JK1]} {\bsc J.M. Johnson and J. Koll\'ar}, {\ninesl K\"ahler-Einstein
metrics on log del Pezzo surfaces in weighted projective 3-space}, Ann. Inst. 
Fourier 51(1) (2001) 69-79.
\item{[JK2]} {\bsc J.M. Johnson and J. Koll\'ar},
{\ninesl Fano hypersurfaces in
weighted projective 4-spaces}, Experimental Math. 10(1) (2001) 151-158.
\item{[KM]} {\bsc J. Koll\'ar, and S. Mori}, {\ninesl Birational Geometry of
Algebraic Varieties}, Cambridge University Press, 1998.
\item{[Mil]} {\bsc J. Milnor}, {\ninesl Singular Points of Complex
Hypersurfaces}, Ann. of Math. Stud. 61, Princeton Univ. Press, 1968.
\item{[MO]} {\bsc J. Milnor and P. Orlik}, {\ninesl Isolated singularities
defined by weighted homogeneous polynomials}, Topology 9 (1970), 385-393.
\item{[Mo]} {\bsc S. Morita}, {\ninesl A Topological Classification of Complex 
Structures on $\scriptstyle{S^1\times \grS^{2n-1}}$}, Topology 14 (1975), 
13-22.
\item{[MU]} {\bsc S. Mukai and H. Umemura}, {\ninesl Minimal Rational 
Threefolds}, LNM 1016, in Algebraic Geometry, M. Raynaud and T. Shioda Eds.,
pgs 490-518, Springer-Verlag, New York, 1983.
\item{[Na]} {\bsc A.M. Nadel}, {\ninesl Multiplier ideal sheaves
and existence of K\"ahler-Einstein metrics of positive scalar curvature}, Ann.
Math. 138 (1990), 549-596.
\item{[Sa]} {\bsc H. Sato}, {\ninesl Remarks Concerning Contact Manifolds}, 
T\^ohoku Math. J. 29 (1977), 577-584.
\item{[TaYu]} {\bsc S. Tachibana and W.N. Yu}, {\ninesl On a Riemannian space
admitting more than one Sasakian structure}, T\^ohoku Math. J. 22
(1970) 536-540.
\item{[Ti1]} {\bsc G. Tian}, {\ninesl K\"ahler-Einstein metrics with positive
scalar curvature}, Invent. Math. 137 (1997), 1-37.
\item{[Ti2]} {\bsc G. Tian}, {\ninesl Canonical Metrics in K\"ahler Geometry}, 
Birkh\"auser, Boston, 2000.
\item{[Us1]} {\bsc I. Ustilovsky} {\ninesl Infinitely Many Contact Structures 
on $\scriptstyle{S^{4m+1}}$}, Int. Math. Res. Notices 14 (1999), 781-791.
\item{[Us2]} {\bsc I. Ustilovsky} {\ninesl Contact Homology and Contact 
Structures on $\scriptstyle{S^{4m+1}}$}, Ph.d. thesis, Stanford Univ., 2000.
\item{[YK]} {\bsc K. Yano and M. Kon}, {\ninesl
Structures on manifolds}, Series in Pure Mathematics 3,
World Scientific Pub. Co., Singapore, 1984.
\medskip
\bigskip \line{ Department of Mathematics and Statistics
\hfil August 2001} \line{ University of New Mexico \hfil }
\line{ Albuquerque, NM 87131 \hfil } \line{ email: cboyer@math.unm.edu,
galicki@math.unm.edu, nakamaye@math.unm.edu\hfil} \line{ web pages:
http://www.math.unm.edu/$\tilde{\phantom{o}}$cboyer,
http://www.math.unm.edu/$\tilde{\phantom{o}}$galicki \hfil}
\bye